%% file: cost.tex
\documentclass[11pt,a4paper]{article}
\usepackage{color}
\usepackage[latin1]{inputenc}
\usepackage{amstext}
\usepackage{graphicx}
\usepackage{amsfonts}
\usepackage{latexsym}
\usepackage{amssymb}
\usepackage{amsthm}
\usepackage{amsmath}
\usepackage{url}

\title{Decentralized Search with Random Costs} 
\author{Oskar
  Sandberg \thanks{The Department of Mathematical Sciences, Chalmers
    University of Technology and Göteborg University. ossa@math.chalmers.se}}

\input{defs.tex}

\begin{document}
\maketitle

\begin{abstract}
  A decentralized search algorithm is a method of routing on a random
  graph that uses only limited, local, information about the
  realization of the graph. In some random graph models it is possible
  to define such algorithms which produce short paths when routing
  from any vertex to any other, while for others it is not.

  We consider random graphs with random costs assigned to the edges.
  In this situation, we use the methods of stochastic dynamic
  programming to create a decentralized search method which attempts
  to minimize the total cost, rather than the number of steps, of each
  path. We show that it succeeds in doing so among all decentralized
  search algorithms which monotonically approach the destination. Our
  algorithm depends on knowing the expected cost of routing from every
  vertex to any other, but we show that this may be calculated
  iteratively, and in practice can be easily estimated from the cost
  of previous routes and compressed into a small routing table. The
  methods applied here can also be applied directly in other
  situations, such as efficient searching in graphs with varying
  vertex degrees.
\end{abstract}

\input{cost_main.tex}

\bibliographystyle{unsrt}
\bibliography{../../tex/ossa}

\end{document}

%% file: defs.tex
\newtheorem{theorem}{Theorem}[section]

\newtheorem{definition}[theorem]{Definition}
\newtheorem{proposition}[theorem]{Proposition}

\newcommand{\PR}{\mathbf{P}}
\newcommand{\given}{\,\mathbf{|}\,}

\newcommand{\E}{\mathbf{E}}

\newcommand{\R}{\mathbb{R}}

%% file: cost_main.tex
\section{Introduction}

Jon Kleinberg introduced the concept of decentralized search
algorithms in his celebrated 2000 paper on the Small-World phenomenon
\cite{kleinberg:smallworld}. In particular, he showed that in certain
random graphs it is possible to find paths between vertices of
poly-logarithmic length even when limited to using only local
knowledge at each step, while in others it is not.

\begin{figure}
   \centering
    \includegraphics[width=12cm]{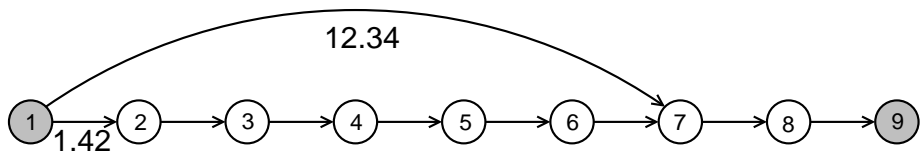} 
    \caption{If the goal is to reach vertex 9 in the fewest steps,
      then certainly vertex 1 should choose the long edge to vertex 7.
      However, if the goal instead is to minimize the total cost,
      routing to 2 might be a better choice.}
\label{fig:examplert}
 \end{figure}

Most of the by now large canon of work in the area (see
\cite{kleinberg:networks} for a recent survey) has been dedicated to
finding and analyzing algorithms that route between two given vertices
in a small number of steps. Typically, the best method in these
situation is greedy: progress to the neighbor which is closest to the
destination. In this paper, we consider a generalized situation where
the cost of passing down an edge is not fixed, but may be a random
variable whose value is known when the choice of where to proceed to
next is made.  The goal is then the minimize the cost of reaching the
destination, which may lead to different priorities when routing (one
may not wish to route to a vertex very close to the destination if the
cost of passing down that edge is very high, see Figure \ref{fig:examplert}).
This is a problem similar to that applied to time process in the field
of stochastic dynamic programming \cite{ross:stochastic}, and we use
similar methods.

The basic idea is this: If a vertex knows the expected cost of each of
his neighbors routing to the destination, and also the cost of him
routing to each of his neighbors, it makes sense for him to choose as
the next step that neighbor which minimizes the sum of these two
costs. While the expected costs are difficult to calculate
analytically, we find that there is little reason to do so. When many
queries are performed, one may start by using any guess as to these
values, and then update these guesses based on past experience. 

We will show analytically that this form of search is well defined,
and that it is optimal among algorithms that monotonically approach
the destination, as well as presenting some results on the order of
the total cost as well as approximation trade-offs. We continue with a
discussion and some experimentation on the practicality of the
approach. Finally, we apply the algorithm with greater generality, and
see that it performs well also in these cases. In particular, we note
1that it can be applied to routing with non-homogeneous degree
distributions.

\subsection{Previous Work}

The original work on decentralized search was done by Kleinberg in
\cite{kleinberg:smallworld} and \cite{kleinberg:navigation}. Much work
has since been done on related problems, in particular further
generalization of the results and improvements on the bounds (see
\cite{barriere:efficient} \cite{kleinberg:dynamics}
\cite{nowell:geographic} \cite{fraigniaud:doubling}
\cite{sandberg:neighbor} for some examples).

In a bid to improve the performance of Freenet
\cite{clarke:protecting}, a decentralized peer-to-peer network, Clarke
\cite{clarke:generation} proposed an algorithm in some ways similar to
that presented here under the name ``NG Routing''. The method was
implemented in Freenet at the time, but has not been used since the
network was re-engineered to route according to the method described
in \cite{sandberg:distributed}. We believe that the technical and
architectural problems experienced by previous versions of Freenet
were unrelated to what is discussed below, and, in light of our
results, that Clarke's ideas were fundamentally sound.

\c{S}im\c{s}ek and Jensen also proposed an algorithm based on the same
principles in \cite{simsek:decentralized}, intended for routing based
on both vertex degree and similarity. Their ``expected-value
navigation'' is based on the same idea as our cost-greedy search,
however they use a rough estimate of the expected routing time, which
cannot be applied to our problem, to make the decisions.  They present
no analytic results. For our take on search in graphs with variable
node degree, see Section \ref{sec:expfurther} below.

For an introduction to the field of stochastic dynamic programming,
see the monograph by Ross \cite{ross:stochastic}. 

\subsubsection{A note on terminology}

The terminology regarding algorithms for decentralized path-finding in
random graphs has not yet settled, and different authors have used
different terms. We have chosen the term ``decentralized search''
(following recent work by Kleinberg \cite{kleinberg:networks}) but
others terms have been used to describe the same thing.
``Navigation'', used in e.g.  \cite{kleinberg:navigation}
\cite{franceschetti:navigation}, is quite common, but perhaps not as
descriptive. We avoid formally calling our algorithms ``routings'' as
this name has previously been used in computer science literature to
describe flow assignments through a graph with limited edge capacities
(for example \cite{stern:decentralized}), a different problem from
that currently studied. However, we do use the terms ``route'' and
``to route'' following their dictionary definitions.

\subsection{Organization}

In Section \ref{sec:def} we set out the basic definitions of 
decentralized search, as well as rigorously define our new
``cost-greedy'' search algorithm. Following this, in Section
\ref{sec:results} we prove the main results regarding cost-greedy
search in networks with sufficient independence. In Sections
\ref{sec:weights} and \ref{sec:practicality} we set out the methods
for applying the results in practice, and, finally, in Section
\ref{sec:exp} we perform simulated experiments to look at the actual
performance of the algorithm.

\section{Definitions}
\label{sec:def}

A cost graph $G = (V, E, C)$ is a graph consisting of a vertex set
$V$, a possibly directed set of edges $E$, and a collection of costs,
$C$.  For each element $(x,y) \in E$, $C(x, y) = C((x,y))$ is an
i.i.d.  positive random variable giving the cost of traveling down
that edge (time taken).

$G$ may be a random graph, by which we mean that for each $x, y \in
V$, there exists a random indicator variable $E_{x,y}$ saying whether
there is an edge from $x$ to $y$. These may be dependent and
differently distributed.

\newcommand{\alga}{\mathcal{A}}

\begin{definition}
  For a given cost graph $G$, a \emph{$z$-search} for a vertex $z \in
  V$, is a mapping $\alga : V \mapsto V$ such that:
  \begin{enumerate}
  \item $\alga(x) = y$ only if $(x, y) \in E$.
  \item $\alga(z) = z$.
  \item For all $x \in V$ there exists is $k < \infty$ such that $\alga^k(x) = z$.
  \end{enumerate}
  A \emph{search} of $G$ is a collection of $z$-routing algorithms
  for all $z \in V$.
\end{definition}

We call $d$ a \emph{distance} on a set $V$ if $d : V \times V \mapsto
\R^+$, if for $x, y \in V$, $d(x, y) = 0$ implies that $x = y$, and if
for $z, y, z \in V$,
\[
d(x, z) \leq d(x, y) + d(y, z).
\]  
A distance is thus a metric without the symmetry requirement. In
particular, any connected digraph $G$ implies a distance $d_G$. For $x
\in V$, $N(x) = \{y \in V : d_G(x,y) = 1\}$ is the set of neighbors of
$x$ in $G$.

\begin{definition}
  A distance $d$ is \emph{adapted for search} in a connected graph $G$
  if for every $x, z \in V$, where $x \neq z$, there exists a $y \in
  N(x)$ such that $d(y,z) < d(x,z)$.
\end{definition}

A distance function $d$ is thus adapted for search if it, in some
sense, reflects the structure of $G$. The most obvious example is of
course $d_G$ itself, but $d$ may also be, for instance, graph distance
on any connected spanning subgraph $H$ of $G$. Another important case
is that if $V$ is a set of points in a metric space, then the space's
metric is adapted for search in $V$'s Delaunay triangulation (see
\cite{sandberg:neighbor}).

\begin{definition}
\label{def:decentralized}
Given a cost graph $G$, and a vertex $z$, a \emph{decentralized
  search} is a $z$-search $\alga$, such that for any $x \in V$ the
random variable $\alga(x)$ is measurable with respect to:
\begin{enumerate}
\item $E_{x,y}$ for all $y \in V$.
\item $C(x,y)$ for all $y \in N(x)$.
\end{enumerate}
\end{definition}

Intuitively, this means that as well as any information about the
graph model, routing at $x$ may use information about which vertices
$x$ does (and does not) have edges to, as well as the costs of passing
down those edges. The definition of decentralized search as originally
given by Kleinberg was slightly broader than this, allowing a route
started at a vertex $x$ to use all the information from $\alga^0(x),
\alga^1(x), \hdots, \alga^{k-1}(x)$ when taking its $k$-th step.
Because our analysis will be restricted to algorithms meeting the
following criteria, excluding this information will make little
difference:

\newcommand{\falga}{\mathcal{F}}

\begin{definition}
  Given a graph $G$, a distance $d$ adapted for search in $G$, and a
  vertex $z$, a \emph{forward search} is a $z$-search $\falga$ such
  that for all $x \in V \backslash \{z\}$
\[
d(\falga(x), z) < d(x, z).
\]
\end{definition}

For a given search $\alga$ and vertex $z$
\[
S^z(x ; \alga) = \inf\left (k : \alga^k(x) = z \right )
\]
is the number of \emph{steps} it takes to reach $z$ from $x$
using $\alga$. Let
\[
T^z(x ; \alga) = \sum_{i = 1}^{S^z(x ; \alga)} C \left
  (\alga^{i-1}(x),\alga^{i}(x) \right )
\]
which is the \emph{cost}.

\subsection{Greedy and Weighted Greedy Search}

\newcommand{\argmin}[1]{\underset{#1}{\text{argmin}}}

In the following, we fix $z$. \emph{Greedy search} is given by
\begin{equation}
\alga_G(x) = \argmin{y \in N(x)}(d(y,z)). 
  \label{eq:gr}
\end{equation}
This is always a decentralized search, and if $d$ is adapted, than it
is a forward search (if $d$ is not adapted, it may not be well defined
to begin with).

Standard greedy search does not take the costs $C$ into account. A
variant that does, is \emph{weighted greedy search}. For a given $z$,
let $w_z(x)$ for $x \in V$ be a collection of weights, where $w_z(z) =
0$.  These weights specify a search algorithm
\begin{equation}
\alga_{w}(x) = \argmin{y \in N(x)}(C(x, y) + w_z(x))
  \label{eq:wgr}  
\end{equation}
which may or may not be well defined. If we restrict this to being a
forward search, we get
\begin{equation}
\falga_{w}(x) = \argmin{y \in N(x)\,:\,d(y,z) < d(x, z)}(C(x, y) + w_z(x))
  \label{eq:fwgr}  
\end{equation}
which is always a well-defined if $d$ is adapted for routing
in $G$.

The behavior of weighted routing depends on $w$. If $w(x)$ is strictly
increasing in $d(x,z)$ then $\alga_w$ and $\falga_w$ are $\alga_G$
when the costs are constant. If $w(x) = 0$ for all $x \in V$, then
$\falga_w$ will simply choose the edge approaching $z$ with the lowest
cost. Since we are trying to minimize the cost of routing, it makes
sense to see $w_z(x)$ as a guess of the cost of reaching $z$ from $x$.
Imagine that we are given a black-box function $f$, where
\[
f(x, \alga) = \E[T^z(x ; \alga)]
\]
that is, $f$ tells us the expected cost of routing from $x$ using any
given algorithm $\alga$. If a vertex $x$ knows that all its neighbors
will use the algorithm $\alga$ to route, it makes sense that it should
want to use the $y$ given by (\ref{eq:wgr}) (or (\ref{eq:fwgr}) if
restricted to forward search) with $w_z(y) = f(y, \alga)$. Extending
this reasoning to any vertex motivates the following definition.

\newcommand{\tw}{{\tilde{w}}}
\begin{definition}
\label{def:tgr}
A forward search $\falga$ is called \emph{cost-greedy search} for $z$
if $\falga = \falga_\tw$, where the weights $\tw$ are given by the
solution to the equations:
\begin{equation}
w(x) = \E[T^z(x ; \falga_w)]\;\;\;{x \in V}
\label{eq:tgr}
\end{equation}
\end{definition}

Equation \ref{eq:tgr} is what is the routing equivalent of what is
known in stochastic dynamic programming texts as the \emph{optimality
  equation}. We will show below that a solution exists in this
context, that this solution is a globally attractive fix-point, and
that cost-greedy search is optimal with respect to expected routing
cost for all forward searches.

\section{Results for Independent Graphs}
\label{sec:results}

Let $d$ be an distance, and $y \leadsto x$ denote the event that there
is an edge from $y$ to $x$ (we allow multiple edges per pair). We let
$G$ be a random graph model on which $d$ is adapted for
routing (with probability 1), and for which $y, x, v, u \in V$ with
$d(y,z) > d(x,z)$ and $d(v,z) > d(u,z)$, $y \leadsto x$ is independent
of $v \leadsto u$ if $y \neq v$. Examples of such graphs are adding
outgoing edges from each vertex with destinations chosen independently
(as in Kleinberg's work \cite{kleinberg:smallworld}) or allowing
each edge, either seen as directed or undirected, to exist
independently of all others (like in classical random graphs and
long-range percolation \cite{newman:percolation}). We call a graph
constructed in this manner \emph{edge independent}.

\begin{theorem}
\label{th:conv}
For an edge independent random graph $G$, there exists a solution
$\tilde{w}$ to equation (\ref{eq:tgr}) so that cost-greedy search is
well defined. $\tilde{w}$ is a globally attractive fix-point of the
iteration given by
\begin{equation}
w_{i+1} =  \E[T^z(x ; \falga_{w_i})].
\label{eq:weights}
\end{equation}
\end{theorem}

\begin{proof}
  Define the rank of $x$ with respect to $z$, $r_z(x)$, as the
  position of $x$ when all the elements in $V \backslash \{z\}$ are
  ordered by increasing distance from $z$, using some deterministic
  tie-breaking rule. Let $r_z(z) = 0$. 

  We will proceed by induction on $r_z(x)$.

  Let $w_0$ be any weighting. Let $x \in V$ such that $r_z(x) = 1$.
  For any forward search $\falga$, $\falga(x) = z$, whence $\E[T^z(x
  ; \falga)]$ does not depend on the algorithm. In particular
  $\E[T^z(x ; \falga_{w})]$ is does not depend on $w$, whence $w_i(x)$
  is constant for $i \geq 1$.

  Let $r_z(x) = k$, and assume that for all $y \in V$ such that
  $r_z(y) < r_z(x)$, $w_i(y)$ takes the same value for all $i \geq
  k-1$. This means that $\E[T^z(x ; \falga_{w_i})]$ takes the same
  value for all $i \geq k$.

  It follows that for all $x \in V$ $w_i(x)$ is fixed for $i \geq
  r_z(x)$. Hence $w_k = w_{k+1}$ for all $k \geq n-1$, and $\tilde{w}
  = w_n$ is a solution to (\ref{eq:tgr}).
\end{proof}

\begin{theorem}
\label{th:optimal}
Let $G$ be an edge independent graph, and $\falga$ a forward search
for a vertex $z \in V$. Then for all $x \in V$
\[
\E[T^z(x ; \falga_\tw)] \leq \E[T^z(x ; \falga)]
\]
where $\falga_\tw$ is cost-greedy search for $z$ as in Definition
\ref{def:tgr}.
\end{theorem}

\begin{proof}
  Use the same definition of $r_z(x)$ as in the proof of Theorem
  \ref{th:conv}. Like there, we will use induction on $r_z(x)$.

  If $r_z(x) \leq 1$, then all forward searches from $x$ are the same,
  and there is nothing to prove. Let $\falga$ be any forward search.
  Given $x \in V$, assume that for all $y \in V$ such that $r_z(y) <
  r_z(x)$, $\E[T^z(y ; \falga_\tw)] \leq \E[T^z(y ; \falga)]$.

  Let $v = \falga(x)$ and $v^* = \falga_T(x)$, the places the
  respective algorithms choose as the next step. Below, we mean by
  ``local knowledge'' that which decentralized algorithm may use, as
  given in \ref{def:decentralized}. We note, crucially, that because
  of our assumptions, $T^z(v ; \falga)$ is independent of local
  knowledge at $x$, while $C(x, y)$ is measurable with respect to it.
  \begin{align*}
\E[ T^z(x ; \falga)] & =  \E[\,\E[T^z(x ; \falga) \given \text{local
  knowledge at }x]\,] \\
    & = \E[\,D(x,v) + \E[T^z(v ; \falga)]\,] \\
& \geq \E[\,D(x,v) + \E[T^z(v ; \falga_\tw)]\,] \\
& \geq \E[\,D(x,v^*) + \E[T^z(v^* ; \falga_\tw)]\,] 
 =  \E[ T^z(x ; \falga_\tw)]
  \end{align*}
  where the first inequality follows by induction since $r_z(v) <
  r_z(x)$ and the last because the expression inside the first
  expectation is what $v^*$ minimizes.
\end{proof}

\subsection{The Small-World Graph}
\label{sec:smallworld}

A particular example graphs meeting the criteria of the last chapter
are the small-world augmentations first introduced by Kleinberg
\cite{kleinberg:smallworld}. This construction starts with a fixed
finite grid $H$, letting $d = d_H$, and creating $G$ by adding a
random outgoing directed edge from each vertex $x$ to destination a
$y$ with probability 
\begin{equation}
\PR(x \leadsto y) \propto 1/d_G(x, y)^\alpha
\label{eq:aug}
\end{equation}
$d$ is naturally adapted for routing in $G$. 

For simplicity, we let $H$ be a ring of $n$ vertices (Kleinberg
originally used a two-dimensional square lattice, but the proofs are
identical). Let $\mathcal{G}(n,\alpha)$ be the family of random graphs
so constructed.

Using previous results about greedy routing on such graphs, we can
calculate the cost order of cost-greedy search. In particular, we can
see the order in $n$ cannot be different from greedy routing.

\begin{theorem}\textbf{(Kleinberg)} If $G \in \mathcal{G}(n,\alpha)$
  with $\alpha = 1$ there exists $N_1$ such that for $n \geq N_1$,
  \[
\E[S^z(x ; \alga_G)] \leq k_1 \log n \log d(x, z)
\]
where $k_1$ is a constant independent of $x$, $z$, and $n$.
\label{th:kleinberg}
\end{theorem}

Further results about such graphs, proved in \cite{barriere:efficient}
and \cite{manku:neighbor} respectively are

\begin{theorem}\textbf{(Barriere et al.)}  If $G \in \mathcal{G}(n,\alpha)$
  with $\alpha = 1$, then there exists $N_2$ such that for $n \geq
  N_2$, 
\[
\E[S^z(x ; \alga_G)] \geq k_2 \log n \log d(x, z)
\]
where $k_2$ is a constant independent of $x$, $z$, and $n$.
\label{th:barriere}
\end{theorem}

\begin{theorem}\textbf{(Singh Manku)} If $G \in
  \mathcal{G}(n,\alpha)$, with $\alpha \geq 0$, then
\[
\E[S^z(x ; \alga_G)] \leq \E[S^z(x ; \alga)]
\]
for any $x, z \in V(G)$ and decentralized search $\alga$.
\label{th:manku}
\end{theorem}

Together, these allow us to prove the observation that

\begin{proposition}
  If $G \in \mathcal{G}(n,\alpha)$ and $0 < \E[C(x,y)] < \infty$ then 
\[
\E[T^z(x ; \falga_{\tilde{w}})] = \Theta(\log n \log d(x, z))
\]
\label{prop:costorder}
\end{proposition}

\begin{proof}
  The upper bound comes directly from Theorems \ref{th:optimal} and
  \ref{th:kleinberg}, since
  \begin{align*}
    \E[T^z(x ; \falga_{\tilde{w}})] & \leq  \E[T^z(x ;
    \alga_G)] \\
& = \E[C(x,y)]  \E[S^z(x ; \alga_G)] \\
& \leq \E[C(x,y)] k_1 \log n \log d(x, y).
  \end{align*}
  where the middle equality follows from the fact $\alga_G$ routes
  independently of the costs, and the simple form of Wald's equation.

  To prove the upper bound, consider all edges in the graph as
  directed, letting the edges of the $H$ be denoted by double directed
  edges. Since a forward search can only ever traverse an edge in one
  direction, this does not affect its cost. Now let
  \newcommand{\Cmin}{C_{\text{min}}}
\begin{equation}
R^z(x ; \falga) = \sum_{i=1}^{S^z(x; \alga)} \Cmin(\falga(x)^{i-1})
\label{eq:rdef}
\end{equation}
where $\Cmin(x) = \min\{C(x,y) : y \in N(x)\}$. This counts, at each
step, the minimum cost of any outgoing edge, rather than the cost of
the edge which was actually used.

Since the degree of each vertex is fixed, $\Cmin(x)$ is i.i.d. for all
$x$. Let $\mathcal{S}_i$ be the $\sigma$-algebra generated by all the
information seen in steps $1,2,\hdots,i$ of the search (as listed in
Definition \ref{def:decentralized}). Note:
\begin{itemize}
\item $\Cmin(\falga^{i+1}(x))$ is independent of $\mathcal{S}_i$.
\item $S_z(x ; \alga)$ is a Stopping Time with respect to the
  filtration $\{S_i\}_{i=1}^\infty$.
\end{itemize}
Thus we may use Wald's Equation to conclude that
\begin{equation}
\E[R^z(x ; \falga)] = \E[\Cmin(x)] \E[S^z(x ; \falga)].
\label{eq:wald}
\end{equation}
We now use the immediate fact that $R^z(x ; \falga) \leq T^z(x ;
\falga)$, followed by (\ref{eq:wald}) and Theorems \ref{th:manku} and
\ref{th:barriere}, to conclude
\begin{align*}
\E[T^z(x ; \falga)] & \geq \E[R^z(x ; \falga)] \\
& =  \E[\Cmin(x)] \E[S^z(x ; \falga)] \\
& \geq \E[\Cmin(x)] \E[S^z(x ; \alga_G)] \\
& \geq \E[\Cmin(x)] k_2 \log n \log d(x,z)
\end{align*}
for sufficiently large $n$. Since this holds for any forward-search
$\falga$, it holds in particular for $\falga_{\tilde{w}}$.
\end{proof}

Proposition \ref{prop:costorder} tells us that in this model, the
order of cost-greedy routing will not be different from that of greedy
routing. The proof of the lower bound assumes, however, that
$\E[C(x,y)] < \infty$ and that the degree of each vertex is bounded as
$n$ grows. Neither of these things, and particularly not the latter,
necessarily hold in applications.

\subsection{Approximated Weights}

We consider the situation when the solution $\tilde{w}$ to
(\ref{eq:weights}) is not known exactly but approximated by another
set of weights.

\begin{proposition}
\label{prop:approx}
  If $\tilde{w}$ is the solution to (\ref{eq:weights}) and $w$ another
  set of positive weights such that
\[
\max_{x \in V} |w(x) - \tilde{w}(x)| \leq \epsilon 
\]
then for any edge-independent graph of size $n$
\begin{equation}
\E[T^z(x ; \falga_w)] - \E[T^z(x ; \falga_{\tilde{w}})] \leq 2 n
\epsilon
\label{eq:nerr}
\end{equation}
and more generally, for any $k \geq 0$
\begin{equation}
\E[T^z(x ; \falga_w)] - \E[T^z(x ; \falga_{\tilde{w}})] \leq 2 \epsilon
\left (k + n \PR(S^z(x ; \falga_w) > k) \right ).
\label{eq:kerr}
\end{equation}
\end{proposition}

\newcommand{\err}{\text{err}_w}

\begin{proof}
  Let $\err(x) = \E[T^z(x ; \falga_w)] - \E[T^z(x ;
  \falga_\tw)]$. It follows that
  \begin{align*}
    \E[T^z(x ; \falga_w)] & = \E[C(x, \falga_w(x)) +
    \E[T^z(\falga_w(x) ;
    \falga_w)]] \\
    & = \E[C(x, \falga_w(x)) + \E[T^z(\falga_w(x) ; \falga_\tw )]] +
    \err(\falga_w(x)) \\
    & =  \E[C(x, \falga_w(x)) + \tw(\falga_w(x))] + \err(\falga_w(x)).\\
    \intertext{Now, since by the definition of a weighted greedy
      search $C(x, F_w(x)) + w(\falga_w(x)) \leq C(x, y) + w(y)$ for
      all $y \in N(x)$} 
    & \leq \E[C(x, \falga_\tw(x)) + w(\falga_\tw(x)) - w(\falga_w(x))
    + \tw(\falga_w(x))]\\
    & \qquad + \err(\falga_w(x)) \\
    & \leq \E[C(x, \falga_\tw) + \tw(\falga_\tw(x))] +
    \E|w(\falga_\tw(x)) - \tw(\falga_\tw(x))| \\
    & \qquad + \E|\tw(\falga_w(x)) - w(\falga_w(x))| + \err(\falga_w(x)) \\
    & \leq \E[T^z(x ; \falga_\tw)] + 2\epsilon + \err(\falga_w(x))
\end{align*}
It follows that for any $k \geq 0$
\begin{equation}
\err(x) \leq 2 \epsilon k + \err(\falga_w^k(x)).
\label{eq:steperr}
\end{equation}
If $k > S^z(x ; \falga_w)$ then $\falga_w^k(x) = z$ and
$\err(\falga_w^k(x)) = 0$, so (\ref{eq:nerr}) follows since $n > S^z(x;
\alga)$ for all searches.

To prove (\ref{eq:kerr}), note that by the same reasoning
\[
\err(\falga^k_w(x)) = \err(\falga_w^k(x) \given \falga_w^k(x) \neq z) \PR(S^z(x ; \falga_w) > k).
\]
Since the graph is edge independent $\err(\falga_w^k(x) \given
\falga_w^k(x) \neq z)$ is simply the error from some point which is
not $z$, but where (\ref{eq:nerr}) still applies.
\end{proof}

What the proposition, and in particular (\ref{eq:steperr}) says is
that if an approximation $w$ of $\tw$ is off by $\epsilon$, then each
step in the routing adds at most $2 \epsilon$ to the optimal routing
time.  This is intuitively clear, since while $\falga_w$ may choose
the wrong vertex in a given step, it can only do so when the total
(actual) cost of routing via that vertex is within $2 \epsilon$ of the
cost of routing via the real one. For the same reason, it is unlikely
that a better bound can be achieved without further assumptions on the
graph and the cost distribution.

\section{Calculating the Weights}
\label{sec:weights}

Theorem \ref{th:conv} provides us with a method of calculating the
weights $\tilde{w}$ for cost-greedy search. One can start by
assigning any initial weighting $w_0$, and then calculate $w_1, w_2,
\hdots$ using (\ref{eq:weights}). 

A closed analytic form for $\E[T^z(x ; \falga_w)]$ as a function of
the vector $w$ is probably very difficult to find, even in the most
simple situations. One can note however that it can be written
recursively as
\[
\E[T^z(x ; \falga_w)] = \sum  \left (
  \E[T^z(x ; \falga_w)] + \E[C(x,y) \given \falga_w(x) = y] \right
)\PR(\falga_w(x) = y)
\]
where the sum is over all $y \in V$ such that $d(y, z) < d(x,z)$. In
the very simplest cases (such as a directed loop with one augmented
outgoing shortcut chord per vertex) it is possible to calculate
$\PR(\falga_w(x) = y)$ and $\E[C(x,y) \given \falga_w(x) = y]$
analytically, in which case $\E[T^z(x ; \falga_w)]$ can be calculated
numerically by recursion. Because this is complicated, and unlikely to
be of much interest in practice, we do not linger on it.

A much more rewarding strategy is to calculate the weights
empirically. That is, start by simulating a large number of searches
from randomly chosen points using $\falga_{w_0}$. While this is being
done, sample the average routing cost to $z$ from each vertex (due to
the Markovian nature of forward search on an edge independent graph, a
vertex may take a sample every time a query passes through it). After
a sufficient number of queries, the average should be an estimate at
$\E[T^z(x ; \falga_{w_0})]$ by the law of large numbers. One may then
take the average costs from each point as $w_1$, and continue in this
manner. Proposition \ref{prop:approx} indicates how close an
approximation is needed, but unfortunately it is not strong enough to
derive a rigorous bound using a polynomial number of samples.

Further, we note two things about the sampling implemented. Firstly,
one needs to be careful about the way the repeated queries are done.
Since we want the edge costs to be random, $C(x,y)$ must be picked
anew, independently, for each query sampled. If the graph is random,
the edges may be redrawn, but must not -- it simply depends on whether
they are to be seen as random or fixed edges in the $G$ above.

Secondly, the proof of Theorem \ref{th:conv} guarantees convergence in
$n$ steps, meaning that an optimal routing is achieved once $w_n$ has
been calculated (if an empirical method is used, the resulting
weighting may still suffer inaccuracies due to the sampling). This is
an unfortunately large number of iterations, especially given that
each may require simulating a large number of queries, but we find
that in practice, much fewer iterations (typically two or three, even
for very large networks) are needed, see Section \ref{sec:exp} below.

\section{Practicality and Decentralization}
\label{sec:practicality}

We proceed to discuss actual applications of Definition \ref{def:tgr}.
On the face of it, the routing method described does not seem
particularly practical. Even if we can calculate the weighting
$\tilde{w}$, this gives a routing table of size of $n$, and such a
table is needed for each $z$ we wish to route for. The complete table
of weights needed to route between any two vertices is thus of size
$n^2$.

Several assumptions can help here however. 

\subsection*{Centralization: Translation Invariance}

If we assume that the graph is translation invariant, then $\E[T^z(x ;
\falga)] = \E[T^0(x - z ; \falga)]$ so $x$ needs only know the routing
cost from each starting vertex to a distinguished vertex $0$. In fact,
in many cases (such as the common case of augmenting single cycle with
random outgoing edges) $\E[T^z(x ; \falga_T)]$ may be exactly, or at
least approximately, a function of $d(x,z)$, in which case $x$ need
only know the expected cost of routing a given distance.  This
knowledge is the same for all $x$, so may be calculated as a single,
global, vector.

\subsection*{Decentralization}

If one wishes for a completely decentralized search system, as, for
instance in peer-to-peer systems such as \cite{clarke:protecting},
then one cannot store a global vector of weights. Instead, each vertex
must store the weights needed to route to every other vertex. In
particular, each vertex $x$ needs to be able to calculate $w_z(y)$ for
each $y \in N(x)$ and $z \in V$. If one assumes translation
invariance, $x$ need only store one such weight vector, and can
translate it to apply to its neighbors. Without such invariance, it
needs to store $|N(x)|$ vectors.

\subsection*{Weight Vector Compression}

In both cases above, however, we are still left with a routing table
size of at least $n$, which is definitely not desirable. The heart of
what makes our method practically useful comes from the fact that the
previous theory about decentralized search makes compression to a
logarithmic size possible.

If we consider graphs of type $\mathcal{G}(n, \alpha)$ described in
Section \ref{sec:smallworld}, we know from Proposition
\ref{prop:costorder} that
\[
E[T^z(x,\falga_{\tilde{w}})] \approx c \E[C] \log(n) \log(d(x,z)).
\]


The utility of this is that if we know that $\E[T^z(x,\falga_T)]$
grows logarithmically with $d(x,z)$ (as indicated by Proposition
\ref{prop:costorder}), we are motivated to assume that it, and thus
the weights $\tw$ in Definition \ref{def:tgr}, can be approximated by
assuming $\tilde{w}(x)$ and $\tilde{w}(y)$ have similar values if $x$
and $y$ are such that $\log(d(x,z)) \approx \log(d(y,z))$. In
particular if $0 \leq d(y,z) - d(x,z) \leq r$ we get
\begin{equation}
| \tilde{w}(y) - \tilde{w}(x) | \approx c_1 \E[C] \log(n) {r \over
  d(x,z)}.
\label{eq:logerr}
\end{equation}
It is easy to prove, using the same methods as in the proof of Theorem
\ref{th:kleinberg} that for greedy routing in $\mathcal{G}(n, \alpha)$
\[
\PR(S^z(x ; \alga_G) \geq \log^3 n) \leq c_2 {\log n \over n}.
\]
Assuming that a similar bound holds for $\falga_\tw$, equation
(\ref{eq:kerr}) with $k = \log^3 n$ in Proposition \ref{prop:approx},
gives that
\begin{align*}
\E[T^z(x ; \falga_w)] - \E[T^z(x ; \falga_{\tilde{w}})] & \leq 4 
\max_{x \in V} |\tw(x) - w(x)| \log^3 n
\end{align*}
so if $r < \epsilon d(x,z) / 4 \log^3 n$ in (\ref{eq:logerr}), then
\begin{align*}
\E[T^z(x ; \falga_w)] - \E[T^z(x ; \falga_{\tilde{w}})] & < \epsilon.
\end{align*}

\begin{figure}[t]
  \begin{center}
    \includegraphics[width=12cm]{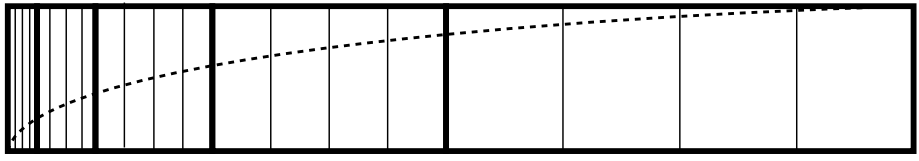} 
  \end{center}
  \caption{Rather than storing every value of $w$ vector (the dotted
    line) as a routing table, we store values at exponentially
    increasing positions, and use these to approximate the values
    between them.}
  \label{fig:rt}    
\end{figure}

Thus, the weight $w(x)$ can be substituted by the weight of a vertex
$\epsilon d(x,z) / \log^3 n$ steps from $x$. To do this, we divide the
routing distances into zones of size $2^i$ for $i = 0,1,2,\hdots$, and
record only the weights of $\log^3 n / \epsilon$ evenly spaced
vertices within each zone (Figure \ref{fig:rt}). The routing table
thus contains a polylogarithmic number of entries ($O(\log^4 n)$), and
yet by using the closest recorded weight as a substitute for $w(x)$,
we incur only an $\epsilon$-error on the total routing cost. Proving
this rigorously, however, depends on tighter bounds then Proposition
\ref{prop:costorder} or even Theorems \ref{th:kleinberg} and
\ref{th:barriere} provide.

We will see experimentally in Section \ref{sec:exp} that a routing
table of size around $E[S^z(x ; \alga_G)]$ works well in practice,
both when using a single vector and in a decentralized system.

While it may seem like a limitation that this will only work on graphs
where routing in a logarithmic number of steps is possible, those are
likely to account for most situations where decentralized search is of
interest. Beyond Kleinberg's small-world model, other cases where
decentralized routing is expected to take a logarithmic number of
steps are hypercubes (where the hamming distance is adapted for
routing), and Chord networks \cite{stoica:chord} (where the circular
distance is).

\section{Experiments}
\label{sec:exp}

\subsection{Direct Applications}

\begin{figure}
  \centering
   \includegraphics[width=7cm]{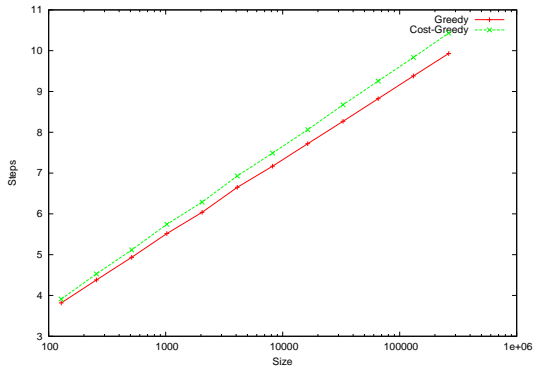} 
   \caption{The steps taken by greedy and cost-greedy search when all
     edge costs are fixed to 1. In theory, the latter should converge
     to the former (which is optimal) but due to the inaccuracy of the
     estimates and the routing-table compression cost-greedy here
     performs about 5\% worse than the optimum.}
  \label{fig:baseline}
\end{figure}

\begin{figure}
  \centering
   \includegraphics[width=7cm]{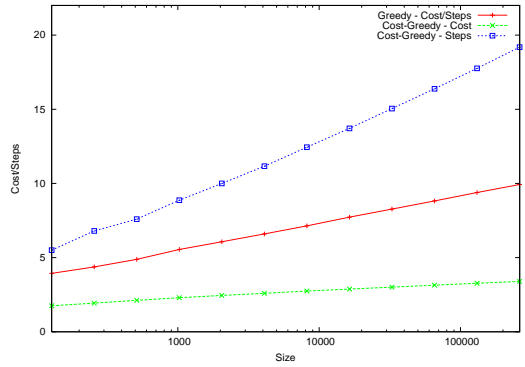} 
   \caption{The cost and steps taken by cost-greedy and greedy search
     on networks with Exp$(1)$ distributed costs along each edge,
     plotted against the size of the network. Networks consist of
     directed rings with $\log_2 n$ directed shortcuts per vertex. }
  \label{fig:bysize}
\end{figure}

 \begin{figure}
   \centering
    \includegraphics[width=7cm]{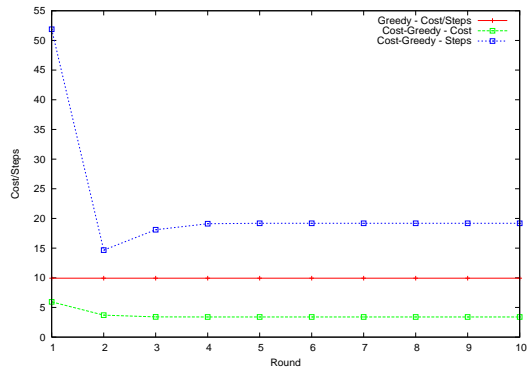} 
    \caption{The cost and steps taken by cost-greedy and greedy
      routing on a network of size 262144 with Exp$(1)$ distributed
      costs along each edge. Cost-greedy performance is plotted
      against iteration of the system in Theorem \ref{th:conv},
      starting with all zero weights.}
   \label{fig:byround}
 \end{figure}

 We start by simulating the algorithm under the most basic conditions.
 We let $G$ consist of a single directed cycle of $n$ vertices,
 augmented with $\log n$ outgoing shortcuts from each vertex,
 according to Kleinberg's small-world model. That is, each shortcut
 from $x$ is to an independently chosen vertex selected with according
 to (\ref{eq:aug}) with $\alpha = 1$, which in this case translates to
 $1 / h_n d(x,y)$, where $d$ is distance in $H$, and $h_n \approx \log
 n$ is a normalizer.

We start by assigning $w_0(x) = 0$ for all $x \in V$, and calculate
the expected routing times by simulating $20 n$ queries between
randomly chosen points (this number of iterations is probably
excessive). We re-sample the costs for each query, but the graph is
kept the same. However, because only one sample vector of the expected
routing times over each distance is kept, we still end up
marginalizing over the shortcuts. We use the logarithmic compression
of $w$ described in Section \ref{sec:practicality} (in practice, we
find this outperforms using a full $w$ vector except when an extremely
large number of queries is simulated), and use $w_{10}$ as an estimate
of the final value.

The difference between cost-greedy and standard greedy search in
terms of query cost depends crucially on the distribution of
$C(x,y)$. Quantitatively, it is possible to make the benefit of cost
greedy as large (or small) as one wishes by a strategic choice of this
distribution. For example, if 
\[
C(x,y) = \begin{cases}
  2 & \text{ with probability } {1 \over 2}\\
  0 & \text{ otherwise.}
\end{cases}
\]
then cost-greedy search will most often incur zero cost assuming the
vertex degree is large enough (as will a simple lowest cost routing).
It would thus be dishonest of us to claim that our methods are
motivated based on the performance achieved with delays chosen by us.
The experiments in this section are thus meant to verify that
cost-greedy search behaves as expected, rather than to illustrate its
benefit: the potential benefit of the algorithm must be evaluated for
every particular situation where it may be applied.

Our first experiment, shown in Figure \ref{fig:baseline}, is thus to
see what happens if we fix the costs to 1 for all the edges. In this
case cost and steps are the same, and since it is known (Theorem
\ref{th:manku}) that greedy search is optimal in the expected number
of steps, the theory tells us that cost-greedy should, ideally, give
the same value.  In fact we find that it under-performs by about 5\%
in all the sizes tested -- presumably due to the empirical estimate of
the expected value, and the losses due to the logarithmic compression
of the weight vector.

Figure \ref{fig:bysize} shows the performance of cost-greedy when the
costs are exponentially distributed as a function of the graph
size. We choose an exponential distribution simply because it is a
common model for waiting times, and the mean of 1 means that the cost
and steps of a route are of the same scale.  We see, as expected, that
cost-greedy search is able to produce routes that cost less by taking
more steps than normal greedy search does. In Figure
\ref{fig:byround} we plot performance for a single network size
against the iterations of (\ref{eq:weights}) when starting with all
zeros. We see that even in a network of hundreds of thousands of
vertices, no measurable performance is gained after the fourth round -
supporting our hypothesis that convergence is a lot faster than the
bound given above.

\subsection{Out-degree Distribution}
\label{sec:expfurther}
 \begin{figure}
   \centering
    \includegraphics[width=7cm]{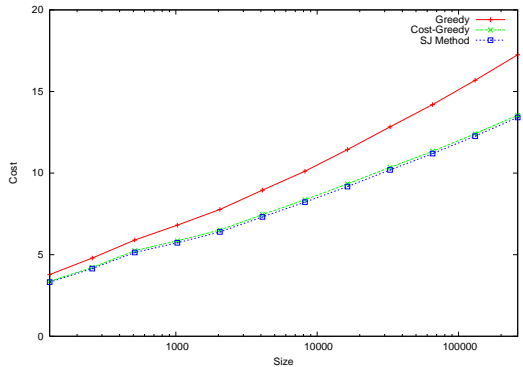} 
    \caption{Example of a network with two different vertex types -
      one tenth of the vertices are augmented by 55 outgoing edges,
      whereas the rest get just 5. All edge costs are exactly one, and
      the values are averaged over four simulations to decrease
      variance. ``SJ Method'' is the method of \c{S}im\c{s}ek and
      Jensen described in Section \ref{sec:expfurther}.}
   \label{fig:deg_dist}
 \end{figure}

 Another question that has been asked about navigability is how to
 route in a network if the vertices have variable degree, and if the
 degree as well as the position of the neighbors is known when the
 routing decision is made. This problem is motivated by the nature of
 social networks, which appear to be navigable, but where it is known
 that the vertex degree follows a heavy-tailed power-law. This problem
 is in many ways similar to that which we discuss above: like with the
 edge costs, degree distributions may incentivize away from a pure greedy
 strategy, and instead call for a trade-off between getting close to
 the destination, and other factors (in this case, wanting to route a
 vertex with high degree).

 \c{S}im\c{s}ek and Jensen \cite{simsek:decentralized} have studied
 this problem by simulation. Their method is fundamentally similar to
 ours: they also seek to choose the neighbor which minimizes the
 expected number of steps to the destination (as we do if the costs
 are fixed to a unit value). However, rather than attempting to
 calculate the fix-point of the weights, they make a rough
 approximation of the value using the right-hand side of the
 inequality
\begin{equation}
\E[S^z(x)] \geq \PR(S^z(x) > 1) = \PR(z \notin N(x))
\label{eq:simsek}
\end{equation}
to estimate the left.

To apply our methods above to the problem, we let the weights be a
function not only of the distance to the destination, but also of the
degree of the vertex. Since this question related only to the number
of steps, we fix all the costs to 1. We expect the weight for any
particular distance to be smaller for vertices with higher degree
(since the amount of ground gained in the first step should be
better).

Figure \ref{fig:deg_dist} shows a simple example of this. In that case
we have exactly two possible out-degrees: a few ($0.1n$) of the
vertices have 55 shortcuts, while the rest have just 5. We compare the
cost-greedy search as used above with regular greedy search and the
method of \c{S}im\c{s}ek and Jensen. The results seems to vindicate
the approximation used in the latter method, with this distribution it
slightly out-performs cost-greedy search, meaning that the numerical
losses in estimating the true weights are greater than the analytic
loss of the approximation. This seems to be the case for most sensible
such distributions, we find that cost greedy search only takes a
slight lead when the popular vertices have more than a hundred times
the degree of unpopular ones. If we presume that cost-greedy search
can come within 5\% of being optimal also here, we are forced to
conclude that so does the SJ method.

One advantage that cost-greedy search has over the method of
\c{S}im\c{s}ek and Jensen, is that their method requires detailed
knowledge about the model in order to calculate the right-hand side of
(\ref{eq:simsek}), which cost-greedy search does not.
Algorithmically, any vertex may implement cost greedy search for its
queries, and it needs only have the ability to to measure the cost of
the queries it sends to its neighbors, nothing more.

\subsection{A Decentralized and Generalized Implementation}
\label{sec:decentralized}
\begin{figure}
   \centering
    \includegraphics[width=7cm]{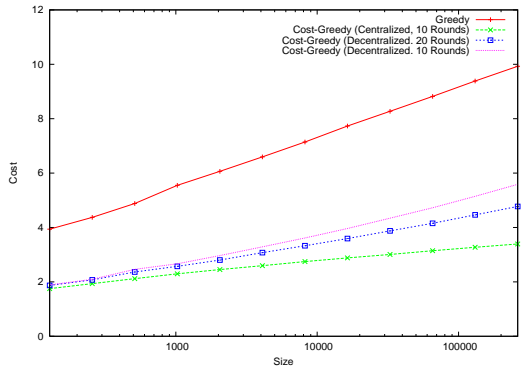} 
    \caption{The same situations as in Figure \ref{fig:bysize}, but
      now including the results when using a separate weight vector at
      every vertex in the manner described in Section
      \ref{sec:decentralized}. A round is $20 \times n$ simulated
      queries as above.}
   \label{fig:decentralized}
 \end{figure}

\begin{figure}
   \centering
    \includegraphics[width=7cm]{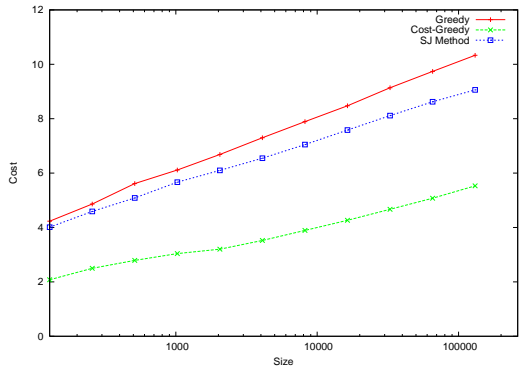} 
    \caption{Using a separate weight vector at every vertex in the
      manner described in Section \ref{sec:decentralized}. Costs are
      Exp$(1)$ distributed, and out-degrees are distributed according
      to a power-law with tail-exponent 2.  ``SJ Method'' is the
      method of \c{S}im\c{s}ek and Jensen described in Section
      \ref{sec:expfurther}.}
   \label{fig:general}
 \end{figure}

 To look at the practical viability of the algorithms described above,
 we also simulate a completely distributed variant. In the
 decentralized variant, we equip each vertex with its own weight
 vector, measuring the mean cost of routing from it to destinations at
 varying distances. Like before we use a $\log_2 n$ compression of the
 weights -- coalescing all distances between $2^k$ and $2^{k+1}$ into
 the same entry -- but unlike above we do not calculate each weight by
 sampling over a fixed number of queries. Instead, we let the weight
 vector at each vertex $x$ be calculated as the mean of the entries in
 a FIFO buffer, showing the cost of the last $m$ queries $x$ has
 routed destined for vertices of every distance category. As before we
 do not change the edges of the graph between the queries from which
 the weights are estimated. Because no marginalizing is occurring
 here, the graph model is actually one \emph{fixed} realization for
 each size -- the expectation is actually taken only over the costs.

 To route a query, $x$ uses the weight vectors of each of the vertices
 in $N(x)$ to minimize (\ref{eq:wgr}). In a real world implementation,
 these values could be periodically copied between neighbors. One
 problem we find with this method is that if $x$ initially has $m$
 queries in a certain distance category that incur a very high cost,
 he will not attract more queries in that category from his neighbors
 (who see it as very costly to send such queries to him), meaning it
 takes a long time to clear the errant values from the FIFO buffer.
 Eventually the buffer will be replaced, if not otherwise then by the
 cost of the queries initiated at $x$ itself, but in our simulations
 we find that this slows the convergence. To alleviate this, we keep a
 count of the number of queries $x$ receives for each distance
 category during an interval. The theory says that these should be
 equal, so if one of the counts has fallen a lot behind (is less than
 a quarter of the queries $x$ receives for itself) we set all the
 values in the buffer to 0.

 Even with this method, the convergence is, as expected, slower than
 in the centralized version. Figure \ref{fig:decentralized} shows the
 equivalent of Figure \ref{fig:bysize} but using local weight vectors
 at each vertex. Here we let $m = 20$. We can note three things: there
 is an absolute performance cost of the decentralized version, the
 cost seems to get worse for larger sizes, but it still considerably
 outperforms greedy search. The first is probably due to
 each estimate of the expected routing time being based on far fewer
 values, while the second is due to us not simulating enough queries
 for full convergence at the large sizes, as seen by the increasing
 difference between the values as 10 and 20 rounds. We note that even
 20 rounds is actually only 400 queries initiated at each vertex -- a
 large number when we must simulate it for a quarter of a million
 vertices, but very little compared to the number of queries one would
 expect in most DHT's or other deployments of distributed networks.

 Finally, in Figure \ref{fig:general} we use the decentralized method
 to route in a situation when we both have exponential edge costs, and
 vertices of varying out-degree (in this case a power-law with
 $\PR(|N(x)| > t) \approx t^{-2}$). Decentralized cost-greedy search
 can optimize both for varying costs and vertex degrees at the same
 time.

\section{Conclusion}

We have presented a method for decentralized search that takes varying
costs of routing down different edges into account. We have showed
that this method is optimal among all such algorithms that
monotonically approach the destination of the query, and that the
necessary weights can be calculated iteratively. On small-world
graphs, we can calculate the order of costs, and say something about
the approximation cost. Beyond these analytic facts, have presented a
number of techniques which make the algorithm practical, and
experimented with actual implementations using simulation.

It would be very desirable to be able to better motivate our
approximations rigorously. To do requires strengthening propositions
\ref{prop:costorder} and \ref{prop:approx}, and perhaps a lot of work
beyond that. In the short term, proving that any polylogarithmic
routing table, and any polynomial number of samples, is sufficient
would be a big improvement.

The long term goal of these studies is to try to find adaptive methods
for decentralized search when the edge costs are not independent
random values, but depend, for instance, on the number of queries that
have passed down the edges recently. Such methods for routing with
congestion are of interest to deployments of peer-to-peer and other
distributed systems.

